\documentclass{article}
\usepackage[a4paper,left=3cm,top=3.5cm]{geometry}
\usepackage{verbatim}
\usepackage[all,cmtip]{xy}

\usepackage{mathtools}

\usepackage{latexsym,bm,easymat}
\usepackage{amsmath,amscd,amssymb,amsfonts,mathdots,ntheorem}
\usepackage{indentfirst}
\usepackage[colorlinks,linkcolor=blue,citecolor=red]{hyperref}
\usepackage{latexsym}
\usepackage{longtable}
\usepackage{cases}
\usepackage{slashed}
\usepackage{tikz-cd}

\usepackage{shorttoc}

\begin{document}
	\makeatletter
	\DeclareRobustCommand\widecheck[1]{{\mathpalette\@widecheck{#1}}}
	\def\@widecheck#1#2{%
		\setbox\z@\hbox{\m@th$#1#2$}%
		\setbox\tw@\hbox{\m@th$#1%
			\widehat{%
				\vrule\@width\z@\@height\ht\z@
				\vrule\@height\z@\@width\wd\z@}$}%
		\dp\tw@-\ht\z@
		\@tempdima\ht\z@ \advance\@tempdima2\ht\tw@ \divide\@tempdima\thr@@
		\setbox\tw@\hbox{%
			\raise\@tempdima\hbox{\scalebox{1}[-1]{\lower\@tempdima\box
					\tw@}}}%
		{\ooalign{\box\tw@ \cr \box\z@}}}
	\makeatother
	
	\def\comp{\ensuremath\mathop{\scalebox{.6}{$\circ$}}}
	\def\QEDclosed{\mbox{\rule[0pt]{1.3ex}{1.3ex}}} %
	\def\QEDopen{{\setlength{\fboxsep}{0pt}\setlength{\fboxrule}{0.2pt}\fbox{\rule[0pt]{0pt}{1.3ex}\rule[0pt]{1.3ex}{0pt}}}}
	\def\QED{\QEDopen} %
	\def\pf{\noindent{\bf Proof}} %
	\def\endpf{\hspace*{\fill}~\QED\par\endtrivlist\unskip \hfill}
	\def\Dirac{\slashed{D}}
	\def\Hom{\mbox{Hom}}
	\def\ch{\mbox{ch}}
	\def\Ch{\mbox{Ch}}
	\def\Tr{\mbox{Tr}}
	\def\End{\mbox{End}}
	\def\Re{\mbox{Re}}
	\def\cs{\mbox{Cs}}
	\def\spinc{\mbox{spin}^\mathbb C}
	\def\spin{\mbox{spin}}
	\def\Ind{\mbox{Ind}}
	\def\Hom{\mbox{Hom}}
	\def\ch{\mbox{ch}}
	\def\cs{\mbox{cs}}
	\def\nTbe{\nabla^{T,\beta,\epsilon}}
	\def\nbe{\nabla^{\beta,\epsilon}}
	\def\mfa{\mathfrak{a}}
	\def\mfb{\mathfrak{b}}
	\def\mfc{\mathfrak{c}}
	\def\mcM{\mathcal{M}}
	\def\mcC{\mathcal{C}}
	\def\mcB{\mathcal{B}}
	\def\mcG{\mathcal{G}}
	\def\mcE{\mathcal{E}}
	\def\mfs{\mathfrak s}
	\def\mcWs{\mathcal W_{\mathfrak s}}
	\def\umcF{\underline{\mathcal F}}
	\def\mce{\mathfrak e}
	\def\mbT{\mathbf T}
	\def\mbc{\mathbf c}
	\def\spinc{$\mbox{spin}^c~$}
	
	\newcommand{\at}{\makeatletter @\makeatother}
	\newtheorem{thm}{Theorem}[section] 
	
	\newtheorem{defi}[thm]{Definition} 
	\newtheorem{lem}[thm]{Lemma} 
	
	\newtheorem{cor}[thm]{Corollary}
	\newtheorem{exa}{Example}[section]
	\newtheorem{prop}[thm]{Proposition}
	\newtheorem{rmk}[thm]{Remark}
	\newtheorem*{ack}{Acknowledgement}
	\newtheorem*{emp}{~~~}
	\newtheorem*{que}{Question}
	\numberwithin{equation}{section}
	\def\slc{\mathfrak{sl}_2(\mathbb C)}	
	\title{Compatible almost complex structures on the Hard Lefschetz condition}
	
	\author{Dexie Lin}
	\date{}
	\maketitle
	\begin{abstract}
		For a compact K\"ahler manifold, it is well-established that its de Rham cohomology satisfies the Hard Lefschetz condition, which is reflected in the equality between the Betti numbers and the Hodge numbers. A special subclass of symplectic manifolds also adheres to this condition. Cirici and Wilson \cite{CW20} employ the variant Hodge number to propose a sufficient criterion for compact almost K\"ahler manifolds to satisfy this condition. In this paper, we show that this condition is only sufficient by presenting examples of compact almost K\"ahler manifolds that fulfill the Hard Lefschetz condition while violating  the equality between the variant Hodge numbers and Betti numbers, that is,
		\[b^1>2h^{1,0}.\]
		This phenomenon contrasts with the behavior observed in compact K\"ahler manifolds.
	\end{abstract}

	
	\\
	
	
	\section{Introduction}
	On a compact  K\"ahler manifold $(M,J,\omega,g)$ with real dimension $2n$,  it is established that  the cup products 
	\begin{equation}
		[\omega^k\wedge]: H^{n-k}(M;\mathbb R)\to H^{n+k}(M;\mathbb R),~k \in \{1, \ldots, n\}\label{eqn-isomorphism}
	\end{equation}  
	are isomorphic, c.f., \cite[Section 6-Chapter 0]{GH}. 
	A  compact symplectic manifold is said to satisfy the Hard Lefschetz condition if the cup products \eqref{eqn-isomorphism} are isomorphisms. This condition serves as a link between K\"ahler manifolds and symplectic manifolds.
	The  isomorphisms in \eqref{eqn-isomorphism} provide a $\slc$-module structure on the de Rham cohomology algebra of $(M,\omega)$. Utilizing this algebraic property of $\slc$-module,  Thurston \cite{Thur76}  constructed a compact symplectic $4$ manifold without K\"ahler structure by violating the Hard Lefschetz condition. 
	This work demonstrated the notable differences between  the topology of symplectic  and K\"ahler manifolds. Such a result disproved Guggenheimer's conjecture \cite{Gug51}. For examples of further symplectic manifolds with non-K\"ahler structure, one can refer \cite{TO97}. 
	
	\noindent
	The Hard Lefschetz condition is also instrumental in the studying of symplectic manifolds.
	\begin{itemize}
		\item  Mathiew \cite{Math95} and Yan \cite{Yan96} showed that the Hard Lefschetz condition is equivalent to establish the symplectic Hodge theory. The concept  of symplectic Hodge decomposition was introduced by Ehresmann and Libermann \cite{EhLib49} and revised  by Brylinski \cite{Bry88}.  Yan's work also disproved Brylinski's conjecture. 
		\item Merkulov \cite[Proposition 1.4]{Mek98} showed that the Hard Lefschetz condition is equivalent to the so called $dd^{\Lambda}$-lemma, see Theorem \ref{thm-ddL-symp-harmonic} for more details. Later, Tomassini and Wang \cite[Theorem B]{TomWang18} generalized this result to complex manifolds with  $J$-symmetric symplectic structures.
		
		\item 	Gompf and Mrowka \cite{GomMro93} constructed symplectic $4$ manifolds satisfying the Hard Lefschetz condition while lacking a  K\"ahler structure. 
	\end{itemize}
	A compact symplectic manifolds is referred to as a Lefschetz manifold, if it satisfies the Hard Lefschetz
	condition. Summarizing the above results, one has the following of inclusions:
	one has the following of inclusions:
	\[
	\mbox{K\"ahler maifolds}\subsetneq \mbox{Lefschetz manifolds}\subsetneq 
	\mbox{Symplectic manifolds}. 
	\]	 This leads to  a natural question in studying compact symplectic manifolds:
	\begin{que}
		What types of compact symplectic manifolds satisfies Hard Lefschetz condition?
	\end{que}
	The algebraic structure of $\slc$-representation gives the natural constraints of the Betti numbers of Lefschetz  manifolds:
	\begin{equation}
		b_{2k+1} \mbox{ are even}, b_{k}\leq b_{k+2},~ b_{2k}>0, \mbox{ for }0\leq k<n-1,
		\label{eqn-Betti-constrains}
	\end{equation}
	Nevertheless, it is possible to construct a compact symplectic $4$ manifold with an even first Betti number that still violates the Hard Lefschetz condition, see Example \ref{exa-non-hlc}.
	Hence, \eqref{eqn-Betti-constrains} is a necessary, but not a sufficient condition, for a symplectic manifold to become a Lefschetz one.
	Cirici and Wilson \cite[Theorem 4.3 and Theorem 5.1]{CW20} provided a sufficient condition, that is  
	\[\sum_{p+q=2 k+1} \ell^{p, q}=2 \sum_{0 \leq p \leq k} \ell^{p, 2 k+1-p} = b^{2 k+1}, \sum_{p+q=2 k} \ell^{p, q}=2 \sum_{0 \leq p<k} \ell^{p, 2 k-p}+\ell^{k, k} = b^{2 k},\]
	\noindent  for compact almost K\"ahler manifolds to admit the Hard Lefschetz condition.  Note that an almost K\"ahler manifold can be regarded as a symplectic manifold equipped with a compatible metric. Here, $\ell^{p,q}$ are variant Hodge numbers, see Theorem \ref{thm-CW-Betti-inequality} in Secrtion 2.

	In this paper, we show that the condition stated above is only sufficient by the following. 
	
	\begin{thm}\label{thm-T4}
		Let $(T^{2n+2},\omega_0)$ be the standard $(2n+2)$-dimensional torus with the standard symplectic form $\omega_0$. Then, for $n>0$ there are infinitely many compatible almost complex structures $J$ such that $$b^1>2  h^{1,0} =2\ell^{1,0}.$$ Moreover, for any integer $0<k\leq n$, there are infinitely many compatible almost complex structures $J$ satisfying
		$$h^{1,0}=\ell^{1,0}=n-k+1.$$
	\end{thm}
	At the same time, Cirici and Wilson {\cite[Theorem 4.1]{CW20}} showed that $\ell^{1,0}=\ell^{0,1}$ is a constraint for compact almost K\"ahler $4$ manifolds.
	We also show that $\ell^{1,0}=\ell^{0,1}$ is only a necessary condition for compact almost Hermitian $4$ manifold becoming almost K\"ahler,  see Theorem \ref{thm-counter}
	
	Additionally, de Bartolomeis and Tomassini \cite{BTom01} gave a sufficient condition in terms of the first non-zero eigenvalues of Hodge Laplacian operators for compact almost K\"ahler manifolds satisfying the Hard Lefschetz condition. 
	For almost K\"ahler $4$ manifiold, we give an analytic condition to hold the equality $b^1=2h^{1,0}$.
	\begin{thm}\label{lemma-b1-almost-kahler}
		Suppose that the inequality
		\[
		\lambda_1>4\|\Delta_{\bar\mu}|_{\mathcal{A}^1}\|,\]
		holds in a compact almost K\"ahler $4$ manifold $(X,J,g,\omega)$,
		where $\lambda_1$ denotes the first non-zero eigenvalue of the Hodge-Laplacian operator $\Delta_d$ restricted on $1$-forms and 
		$\|\Delta_{\bar\mu}|_{\mathcal{A}^1}\|$ denotes the operator norm of $\Delta_{\bar \mu}:L^2(\mathcal{A}^1)\to L^2(\mathcal{A}^1)$.
		Then,  $b^1 =2\ell^{1,0}=2h^{1,0}$.
	\end{thm}      
Note that for any compact almost K\"ahler manifold, the Hodge number  $ {h}^{1,0}$ is equal to $\ell^{1,0}$
the generalized Hodge number  defined  by Cirici and Wilson  \cite{CW20}, i.e., $ {h}^{1,0}=\ell^{1,0}$ by Lemma \ref{lemma-h-ell}. For their definitions, refer to Section 2. 



The structure of this paper is as follows: In Section 2.1, we review the basic notions and properties of the Hard Lefschetz condition. In Section 2.2, we review the  Dolbeault cohomology on almost complex manifold and almost K\"ahler manifolds and the related properties. In  Section 3,  we give the construction of the almost complex structure in Theorem \ref{thm-T4}; ; in the end we give proof of Theorem \ref{lemma-b1-almost-kahler}.

\begin{ack}
	The author would thank Teng Huang and Hongyu Wang for introducing this topic and valuable discussion, and Adriano Tomassini for suggesting considering Example \ref{exa-non-hlc}. The author also thanks Yi Wang for pointing the mistakes in the previous version. This work is partially supported by  National Natural Science Foundation of China No. 12301061 and Fundamental Research Fundsfor the Central Universities No. 2023CDJXY-043. 
\end{ack}

\section{ Hard Lefschetz condition and Almost complex manifolds}
In this section, we review some necessary and related notions and definitions on the Hard Lefschetz condition. 
\subsection{Hard Lefschetz condition on linear spaces}
In this subsection, we give a quick review about the Hard Lefschetz condition on linear symplectic vector space and some properties on the toy model of finite dimensional spaces, c.f., \cite[Section 2 and 3]{TomWang18}. 
First, we review some basic notions about the  symplectic space. 
\begin{defi}
	Let $V$ be a $2n$-dimensional real vector space. A bilinear form $\omega$ on $V$ is called symplectic, if  $\omega$ is skew-symmetric and non-degenerate. 
\end{defi}
Similar to the metric case, 
the symplectic form $\omega$ defines a star operator $*_s: \wedge^p V^* \rightarrow \wedge^{2 n-p} V^*$   by
$$
\mu \wedge *_s \nu=\omega^{-1}(\mu, \nu) \frac{\omega^n}{n !},
$$
where $\omega^{-1}$ is the associated symplectic form on $V^*$.
Set $L:\Lambda V^*\to \Lambda V^*$ as $\alpha\mapsto \omega \wedge\alpha$ and $\Lambda:=*_s^{-1} L *_s$. 



\noindent
It is not hard to see, 
$$
\omega^{-1}(L u, v)=\omega^{-1}(u, \Lambda v),
$$
i.e., $\Lambda$ is the adjoint of $L$.   For the convenience of the later arguments,  we set  $L^r:=\frac{\omega^r}{r !}$.


\begin{defi}
	A linear map $J: V \rightarrow V$ is called   an almost complex structure on $V$ if $J^2=-Id$.
	An almost complex structure $J$ is said to be tamed by $\omega$ if
	$$
	\omega(u, J u)>0,
	$$
	for each non-zero $u \in V $. 
	We say $J$ is $\omega$-compatible if it is both taming and symmetric, i.e., $\omega(u, J v)=\omega(v, J u),
	$
	for each $u, v \in V$.
	
\end{defi} 


An almost complex structure $J$ induces the decomposition 
\[V_{\mathbb C}:=V\otimes{\mathbb C}=V^{1,0}\oplus V^{0,1},\]
where $V^{1,0}$ and $V^{0,1}$ are the eigenspaces of $J$ corresponding to eigenvalues $i$ and $-i$ respectively.
This $J$  action  can be extended to any $p$-form $\psi$   by 
$$(J\psi)\left(v_1, \cdots, v_p\right)= \psi\left(J v_1, \cdots, J v_p\right).$$
Similarly, we  decompose $V^*_{\mathbb C}:=V^{*}M\otimes_{\mathbb R}\mathbb C$ as
$V^*_{\mathbb C}=V^{*,(1,0)}_{\mathbb C}\oplus V^{*,(0,1)}_{\mathbb C}$.
Let  $\bigwedge^{p, q} V^*$ denote the space  of   forms of $\bigwedge^{p,q}V^*:=\bigwedge^pV^{*,(1,0)}_{\mathbb C}\otimes\bigwedge^qV^{*,(0,1)}_{\mathbb C}$. 
An  almost complex structure $J$ on $V$  compatible with $\omega$ 
defines an Hermitian inner product structure on $\bigwedge^{p, q} V^*, 0 \leq p, q \leq n$.
The associated Hodge star operator $*: \bigwedge^{p, q} V^* \rightarrow \wedge^{n-q, n-p} V^*$ is defined by
$$
u \wedge * \bar{v}=(u, v) \omega^n/n! .$$
The following identity gives a relation between $*$ and $*_s$,  $$
*=*_s \comp J=J \comp*_s .
$$
The above relation shows that any $*$-harmonic form is also $*_s$-harmonic form if and only if the compatible structure $J$ preserves the space of $*$-harmonic forms. 
Now we give the definition of Lefschetz space and $\slc$-module, more details can be found in \cite[Section 6-Chapter 0]{GH}. 
\begin{defi}		Let $A=\bigoplus_{k=0}^{2 n} A^k$ be a direct sum of complex vector spaces with possibly infinite dimensions. We say that $L \in \operatorname{End}(A)$ satisfies the Hard Lefschetz Condition and $(A, L)$ is a Lefschetz space if		$$		L\left(A^l\right) \subset A^{l+2},~ L\left(A^{2 n-1}\right)=L\left(A^{2 n}\right)=0,		$$		for $ 0 \leq l \leq 2(n-1)$		and each $L^k: A^{n-k} \rightarrow A^{n+k}, 0 \leq k \leq n$   is an isomorphism.		
\end{defi} 




We review the definitions of Lefschetz space and $\slc$-triple. 
\begin{defi}
	Set $\Lambda=*_s^{-1} L *_*, B:=[L, \Lambda]$.   $(L, \Lambda, B)$ is called an  $\slc$-triple on $(A, L)$.
\end{defi}
The general K\"ahler identity implies the following result.

\begin{thm}[c.f.~{\cite{Gui01}}]
	Let $(A, L)$ be a Lefschetz space. Let $d$ be a $\mathbb C$-linear endomorphism of $A$ such that $d\left(A^l\right) \subset A^{l+1}$ and $d^\Lambda:=(-1)^{l+1}*_sd*_s:A^l\to A^{l-1}$. If $[d, L]=0$, then
	\begin{itemize}
		\item[(1)] $( \ker\left.d \cap \operatorname{ker} d^{\Lambda}, L\right)$ and  $(\ker\left.d \cap \operatorname{ker} d^{\Lambda}, \Lambda\right)$ are Lefschetz.
		\item[(2)]  $\left(\operatorname{im} d+\operatorname{im} d^{\Lambda}, L\right)$ and $\left(\operatorname{im} d+\operatorname{im} d^{\Lambda}, \Lambda\right)$ are Lefschetz.
		\item[(3)]  Assume further that $d^2=0$. Then, $\left(\operatorname{im} d d^{\Lambda}, L\right)$ and $\left(\operatorname{im} d d^{\Lambda}, \Lambda\right)$ are Lefschetz.
	\end{itemize}

\end{thm}

Let $(A, L)$ be a Lefschetz space. Let $d$ be a $\mathbb C$-linear endomorphism of $A$ such that $d\left(A^l\right) \subset A^{l+1}$. We call $(A, L, d)$ a Lefschetz complex if $d^2=0$.
Set 

$$
H_d=\bigoplus_{k=0}^{2 n} H_d^k, \quad H_d^k:=\frac{\operatorname{ker} d \cap A^k}{\operatorname{Im} d \cap A^k},
\quad 
H_{d^{\Lambda}}=\bigoplus_{k=0}^{2 n} H_{d^{\Lambda}}^k, \quad H_{d^{\Lambda}}^k:=\frac{\operatorname{ker} d^{\Lambda} \cap A^k}{\operatorname{Im} d^{\Lambda} \cap A^k} .
$$

One has the following theorem. 


\begin{thm}[c.f. {\cite[Theorem 3.5]{TomWang18}}]\label{thm-ddL-symp-harmonic}
Let $(A, L, d)$ be a Lefschetz complex. Assume that $[d, L]=0$. Then the followings are equivalent:
\begin{itemize}
\item[(1)]  $(A, L, d)$ satisfies the $d d^{\Lambda}$-Lemma.
\item[(2)] The natural map  $\ker d \cap\ker d^{\Lambda} \rightarrow H_d$ is surjective and the  $\slc$-triple $(L, \Lambda, B)$ on  $\ker d \cap \ker d^{\Lambda}$ induces the $\slc$-triple on $\mathrm{H}_d$.
\item[(3)] $\left(H_d, L\right)$ satisfies the Hand Lefschetz condition.
\item[(4)] $\left(H_{d^{\Lambda}}, \Lambda\right)$ satisfies the Hard Lefschetz condition.
\end{itemize}

\end{thm} 
Here a Lefschetz complex  $(A, L, d)$ satisfies the $d d^{\Lambda}$-Lemma if
$$
\ker d \cap \ker d^{\Lambda} \cap\left(\operatorname{im} d+\operatorname{im} d^{\Lambda}\right)=\operatorname{im} d d^{\Lambda},
$$
on each $A^k, 0 \leq k \leq 2 n$.

For a compact symplectic manifold $(M,\omega)$, setting $A=\mathcal{A}^*$, $L=\omega\wedge$ and $d$=de Rahm differential, we say $(M,\omega)$ satisfies the Hard Lefschetz condition or $(M,\omega)$ is a Lefschetz manifold, if it satisfies one of the equivalent statements in Theorem \ref{thm-ddL-symp-harmonic}. 

\subsection{Almost complex structures}
In this subsection, we review the  Dolbeault cohomology on almost complex manifolds and generalized Hodge number on almost K\"ahler manifolds. We first review some basic definitions and notions. 
\begin{defi}
For a smooth manifold $M$ of real dimension $2n$. If there exists  a smooth section $J\in\Gamma(End(TM))$ on $TM$  satisfying $J^2=-1$, then we call
$(M,J)$ an almost complex manifold with almost complex structure $J$.
\end{defi}
Similar to the decompositions on linear toy model $V$ of the previous subsection,  the  parallel decompositions on
$T_{\mathbb C}M=TM\otimes_{\mathbb R}\mathbb C$ and $T^*_{\mathbb C}M:=T^{*}M\otimes_{\mathbb R}\mathbb C$ also hold. 
Denoting $ \mathcal A^*(M)$ by the space of sections of  $\bigwedge^*T^*_{\mathbb C}M$,   it gives that:
\[ \mathcal A^r(M)=\bigoplus_{p+q=r} \mathcal A^{p,q}(M),\]
where $\mathcal{A}^{p,q}(M)$ denotes the space  of   forms of $\bigwedge^{p,q}T^*_{\mathbb C}M:=\bigwedge^pT^{*,(1,0)}_{\mathbb C}M\otimes\bigwedge^qT^{*,(0,1)}_{\mathbb C}M$.
By 
the formula \[d\psi\in   \mathcal{A}^{p+2,q-1}(M)+   \mathcal{A}^{p+1,q}(M)+
\mathcal{A}^{p,q+1}(M)+   \mathcal{A}^{p-1,q+2}(M),\]
for any $\psi\in   \mathcal{A}^{p,q}(M)$,
the de Rham differential can be decomposed  as follows:
\[
d=\mu+\partial+\bar{\mu}+\bar{\partial} \mbox{ on } \mathcal{A}^{p,q},
\]
where   each  component is a derivation. The bi-degrees of the four components are  given by
\[|\mu|=(2,-1)
,~|\partial|=(1,0),~|\bar{\partial}|=(0,1),  ~|\bar{\mu}|=(-1,2).
\]
Observe that the operators $\partial,~\bar\partial$ are of the first order and the operators $\mu,~\bar\mu$   are of the zero order. 
Similar to the case of Riemannian manifolds, we introduce the definitions of almost Hermitian manifold and almost K\"ahler manifold.
\begin{defi}
A metric $g$   on an almost complex manifold $(M,J)$ is called    Hermitian, if $g$ is $J$-invariant, i.e., $g(J-,J-)=g(-,-)$.
The imaginary component $\omega$ of $g$ defined by $\omega=g(J-,-)$, is a real non-degenerate $(1,1)$-form. 
The quadruple   $(M,J,g,\omega)$ is called an \textit{almost Hermitian} manifold. 
Moreover, when $\omega$ is $d$-closed, the quadruple   $(M,J,g,\omega)$ is called  an
\textit{almost  K\"ahler} manifold. 
\end{defi}
Recall that by doing polar decomposition, one can equip an almost K\"ahler structure on any symplectic manifold.
The   Hermitian metric $g$ induces  a unique $\mathbb{C}$-linear Hodge-$*$ operator 
$$
*: \Lambda^{p, q} \rightarrow \Lambda^{n-q, n-p},
$$
$\mbox{defined by }
g\left(\phi_1, \phi_2\right) d V=\phi_1 \wedge *\overline{ \phi_2}$, 
where $dV$ is the volume form $\frac{\omega^n}{n!}$ and $\phi_1, \phi_2 \in \Lambda^{p, q}$. For any two forms $\varphi_1\mbox{ and }\varphi_2$,  we define their  product  by, 
$$(\varphi_1,\varphi_2)=\int_M\varphi_1\wedge *\overline{\varphi_2}.$$  
Let  $\delta^*$ denote the formal-$L^2$-adjoint of $\delta$ for $\delta\in\{\mu,\partial,\bar\partial,\bar\mu,d\}$ with respect to the metric $g$. 
By the integration by part, one has the identities 
\[\delta^*=-*\bar\delta*  \mbox{ for }\delta\in\{\mu,\partial,\bar\partial,\bar\mu,d\}.\] 
By the identity $\Delta_{\bar\partial}+\Delta_\mu=\Delta_\partial+\Delta_{\bar\mu}$ on almost K\"ahler manifold, c.f. \cite[Lemma 3.6]{BTom01}, one has  the following theorem.  



\begin{thm}[Cirici and Wilson {\cite[Theorem 4.1]{CW20}}]\label{thm-CW-02}
On any compact almost K\"ahler manifold with dimension $2 n$, it holds
$$
\mathcal{H}_d^{p, q}=\mathcal{H}_{\bar\partial}^{p, q} \cap \mathcal{H}_\mu^{p, q}=\mathcal{H}_{\partial}^{p, q} \cap \mathcal{H}_{\bar{\mu}}^{p, q}=\ker(\Delta_{\bar\partial}+\Delta_\mu),
$$
for all $(p, q)$. Moreover,
we also have the following identities:
\begin{itemize}
\item[$(1)$](Complex conjugation)
$$
\mathcal{H}_{\bar{\partial}}^{p, q} \cap \mathcal{H}_\mu^{p, q}=\mathcal{H}_{\bar{\partial}}^{q, p} \cap \mathcal{H}_{{\mu}}^{q, p} .
$$
\item[$(2)$]  (Hodge duality) 
$$
*: \mathcal{H}_{\bar{\partial}}^{p, q} \cap \mathcal{H}_\mu^{p, q} \rightarrow \mathcal{H}_{\bar{\partial}}^{n-q, n-p} \cap \mathcal{H}_\mu^{n-q, n-p} .
$$
\item[$(3)$] (Serre duality)
$$
\mathcal{H}_{\bar\partial}^{p, q} \cap \mathcal{H}_\mu^{p, q} \cong \mathcal{H}_{\bar\partial}^{n-p, n-q} \cap \mathcal{H}_\mu^{n-p, n-q} .
$$
\end{itemize}
\end{thm}

For the convenience of  later arguments in the next section, we set $\mathcal{H}^{p,q}_{\bar\partial,\mu}=\mathcal{H}^{p,q}_{\bar\partial }\cap \mathcal{H}^{p,q}_{\mu}$ and $\ell^{p,q}=\dim_{\mathbb C} \mathcal{H}^{p,q}_{\bar\partial,\mu}$.
On the other hand, expanding  the relation $d^2=0$  implies  the following identities:
\[
\bar{\mu} \bar{\partial}+\bar{\partial} \bar{\mu}=0, ~
\bar{\mu}^{2}=0, \mbox{ and }\bar \partial^2|_{\ker(\bar\mu)}\equiv0(\bmod ~ im(\bar\mu)). \]
These identities induce a  Dolbeault-type cohomology. 
\begin{defi}[Cirici and Wilson {\cite[Definition 3.1]{CW21}}]
The Dolbeault cohomology of an almost  complex $2n$-dimensional manifold
$(M,J)$ is given  by
\[
H^{p,q}_{Dol}=H^q(H^{p,*}_{\bar\mu},\bar\partial)=
\frac{\ker(\bar\partial:H^{p,q}_{\bar\mu}\to H^{p,q+1}_{\bar\mu})}
{im(\bar\partial:H^{p,q-1}_{\bar\mu}\to H^{p,q}_{\bar\mu})},
\]
where $H^{p,q}_{\bar\mu}=\frac{\ker(\bar\mu:\mathcal{A}^{p,q}\to \mathcal{A}^{p-1,q+2})}{im(\bar\mu:\mathcal{A}^{p+1,q-2}\to \mathcal{A}^{p,q})}$.
\end{defi}
Note that
Cirici and Wilson also gave the definition in terms of the Hodge-type filtration.
\begin{thm}[Cirici and Wilson {\cite[Theorem 4.3 and Theorem 5.1]{CW20}}]\label{thm-CW-Betti-inequality}
For any compact almost K\"ahler manifold of dimension $2 n$, the following are satisfied:
\begin{itemize}
\item[(1)]  In odd degrees, it establishes
$$
\sum_{p+q=2 k+1} \ell^{p, q}=2 \sum_{0 \leq p \leq k} \ell^{p, 2 k+1-p} \leq b^{2 k+1} .
$$
\item[(2)]   In even degrees, it holds that
$$
\sum_{p+q=2 k} \ell^{p, q}=2 \sum_{0 \leq p<k} \ell^{p, 2 k-p}+\ell^{k, k} \leq b^{2 k},
$$
with $\ell^{k, k} \geq 1$ for all $k \leq n$.
\end{itemize}
Moreover, if the equalities $2 \sum_{0 \leq p<k} \ell^{p, 2 k-p}+\ell^{k, k}=b^{2 k}$ and $2 \sum_{0 \leq p \leq k} \ell^{p, 2 k+1-p} =b^{2k+1}$ hold, then this almost K\"ahler manifold satisfies the Hard Lefschetz condition.
\end{thm}	
We end this section by showing that the Hodge number $h^{1,0}$ is equal to  $\ell^{1,0}$ on any compact almost K\"ahler manifold.

\begin{lem}\label{lemma-h-ell}
Let $(M,g,J,\omega)$ be a compact almost K\"ahler manifold. Then, it holds that 
$$h^{1,0} =\ell^{1,0}. $$
\end{lem}
\begin{pf}
By the identities $\Delta_{\bar{\partial}}+\Delta_\mu=\Delta_\partial+\Delta_{\bar\mu}$ and   $\Delta_\mu\mid_{\mathcal{A}^{1,0}}=0$,   any $u\in {H}^{1,0}_{Dol}$ implies $u\in\ker(\Delta_\partial+\Delta_{\bar\mu})=\ker(\Delta_\partial)\cap\ker(\Delta_{\bar\mu})$. 

Conversely, any $u\in\mathcal{H}^{1,0}_{\bar\partial,\mu}=\mathcal{H}^{1,0}_{\partial,\bar\mu}$  induces $u\in\ker(\bar\partial)\cap\ker(\bar\mu)$ by the same arguments.
\end{pf}

\section{ The proof}
First, let us warm up on a simple case  $T^4$, i.e. $n=2$.
\noindent
Consider $(T^4,\omega_0)$ as the standard $4$-dimensional torus  with the standard symplectic form $\omega_0$. Let $(x^1,x^2,x^3,x^4)$ be the standard coordinate of $T^4$. Fix a positive smooth function $q$ depending only on $(x^3,x^4)$.  Define the almost complex structure $J$ on $T^4$ by
\[J\frac{\partial}{\partial x^1}=\frac{1}{q} \frac{\partial}{\partial x^2},~J\frac{\partial}{\partial x^2}=-q\frac{\partial}{\partial x^1},~J\frac{\partial}{\partial x^3}=\frac{\partial}{\partial x^4},~J\frac{\partial}{\partial x^4}=-\frac{\partial}{\partial x^3}.\]
It is clear that the standard symplectic $\omega_0=dx^1\wedge dx^2+dx^3\wedge dx^4$ of $T^4$ is compatible with $J$.
Consider the   $(1,0)$-forms 
$$\theta^1=dx^1+iqdx^2,~\theta^2=dx^3+idx^4.$$
The dual of $\{\theta^1,\theta^2\}$ is 
\[\{v_1=\frac{1}{2}(\frac{\partial}{\partial x^1}-\frac{i}{q}\frac{\partial}{\partial x^2}),v_2=\frac{\partial}{\partial z^2}=\frac{1}{2}(\frac{\partial}{\partial x^3}-i\frac{\partial}{\partial x^4})\}.\]

Using the structure equations, we obtain
\begin{eqnarray}
\bar\mu(\theta^1) 
=\frac{1}{2q}\bar v_2 q\bar\theta^1 \wedge \bar\theta^2,\label{eqn-T4-1}\\
\bar\partial(\theta^1)=-\frac{1}{2q}\frac{\partial q}{\partial z^2}\theta^2\wedge\bar\theta^1-\frac{1}{2q}\frac{\partial q}{\partial \bar z^2}\theta^1\wedge\bar\theta^2,\label{eqn-T4-2}
\end{eqnarray}
and
\begin{eqnarray}
\partial(\theta^1)=-\frac{1}{2q}v_2q\theta^1\wedge\theta^2.\label{eqn-T4-3}
\end{eqnarray}


Write
\[T^{2n+2}=(T^2)^{\times k}\times T^{2n+2-2k}.\]
Now, we give the proof by perturbing the almost complex structure on the first $2k$-space. 

\begin{pf}\textbf{(of Theorem \ref{thm-T4})}

It suffices to show the second statement. Let $k$ be a fixed integer. We define the almost complex by:
\[J\frac{\partial}{\partial x^{2i+1}}=\frac{1}{q}\frac{\partial}{\partial x^{2i+2}},~J\frac{\partial}{\partial x^{2i+2}}=-q\frac{\partial}{\partial x^{2i+1}},\]
for $0\leq i<k$, where $q$ is a positive smooth function depending only on $(x^{2n+1},x^{2n+2})$. For $k\leq j\leq n$, we set
\[J\frac{\partial}{\partial x^{2j+1}}=\frac{\partial}{\partial x^{2j+2}},~J\frac{\partial}{\partial x^{2j+2}}=-\frac{\partial}{\partial x^{2j+1}}.\]
Let the critical set of $q$ be of nowhere dense, for example, consider $q=e^{\sin(x^{2n+1}+x^{2n+2})}$. 
It is not hard to observe that $J$ is compactible with the standard symplectic form
\[\omega_0=dx^1\wedge dx^2+...+dx^{2n+1}\wedge dx^{2n+2}.\]
Next, we define the co-frame 
$\{\theta^1,...,\theta^k,\theta^{\alpha}\}_{k+1\leq\alpha\leq n+1}$ as follows:
\[\theta^j=dx^{2j-1}+iq dx^{2j},\]
for $1\leq j\leq k$ and 
\[\theta^\alpha=dx^{2\alpha-1}+idx^{2\alpha},\]
for $k<\alpha\leq n+1$.
The dual vector fields are given by:
\[v_j=\frac{1}{2}(\frac{\partial}{\partial x^{2j-1}}-\frac{i}{q}\frac{\partial}{\partial x^{2j}}),\]
for $1\leq j\leq k$, and 
\[v_\alpha=\frac{1}{2}(\frac{\partial}{\partial x^{2\alpha-1}}-i\frac{\partial}{\partial x^{2\alpha}}),\]
for $k<\alpha\leq n+1$.
Analogous to \eqref{eqn-T4-1} and \eqref{eqn-T4-2} on $T^4$, we obtain the structure equations 
\begin{eqnarray}
\bar\mu(\theta^{j+1})=\frac{\bar v_{n+1}q}{2q}\bar\theta^{j+1}\wedge \bar\theta^{n+1}, \label{eqn-Tn-str-1}\\
\bar\partial(\theta^{j+1})=-\frac{v_{n+1}q}{2q} \theta^{n+1}\wedge\bar\theta^{j+1}-\frac{\bar v_{n+1}q}{2q} \theta^{j+1}\wedge\bar\theta^{n+1},
\end{eqnarray}
for $0\leq j\leq k-1$. Let $u=\sum_{l=1}^{n+1}f_l\theta^l$ be a smooth $(1,0)$-form in $\ker({\bar\partial})\cap\ker({ \bar\mu})$. 
From $\bar\mu(\sum_lf_l\theta^l)=0$ and \eqref{eqn-Tn-str-1}, we derive
\[\bar v_{n+1}q\cdot f_j=0,~1\leq j\leq k.\]
Since $\bar v_{n+1}q$ is a.e.  non-zero and $f_j$ is smooth, $f_j\equiv0$  for $1\leq j\leq k$. That is,
$$u=\sum_{k< l\leq n+1}f_l\theta^l.$$ 
To analyze the equation $\bar\partial u=0$, we expand it as
\begin{eqnarray}
0= \sum_{k<\alpha\leq n+1} (\sum_{1\leq j\leq k}\bar v_j f_\alpha \bar\theta^j\wedge\theta^\alpha)+\bar\partial_{n+1-k}(\sum_{k<\alpha \leq n+1} f_\alpha\theta^\alpha),\label{eqn-Tn-dbar}
\end{eqnarray}
where $\bar\partial_{n+1-k}$ is the standard $\bar\partial$-operator of the  K\"ahlerian   $T^{2n+2-2k}$  with respect to the standard complex structure. By employing a Fourier expansion for $f_\alpha$(see \cite[Theorem 3.57]{EW17}), express it as:
\begin{equation}
f_\alpha=\sum_{n_1,...,n_{2k}}a_{n_1,...,n_{2k};\alpha}(x^{2k+1},...,x^{2n+2})e^{2\pi i n_1\cdot x^1}\cdots e^{2\pi in_{2k}\cdot x^{2k}}.  \label{Fourier-1}
\end{equation}
From the first term on the right-hand side of \eqref{eqn-Tn-dbar}, we conclude that
\begin{equation}
\bar v_{j}f_\alpha=0,\label{dbar}
\end{equation}
for $1\leq j\leq k$.
Substituting \eqref{Fourier-1} into \eqref{dbar},  we obtain
\[ \pi i\cdot (n_{2j-1}+\frac{i}{q}n_{2j})\cdot a_{n_1,...n_{2k};\alpha}(x^{2k+1},...,x^{2n+2})=0.\]
Thus, $n_j=0$ for each $j\in\{1,...,2k\}$, i.e. $f_\alpha$  depends only on $(x^{2k+1},...,x^{2n+2})$. 
Consequently, we reduce our problem to solving the standard $\bar\partial$-equations  on $T^{2n+2-2k}$. Therefore, $f_\alpha$'s are all constants and $\ell^{1,0}=h^{1,0}=n+1-k$ by Lemma \ref{lemma-h-ell}. 

To finish the proof, it is clear that
there are infinitely many choices of $q$ with nowhere dense critical sets, such as the Morse functions. This means there are also infinitely many compatible almost complex structures that satisfy  the desired property. 
\end{pf}

By the similar idea, we have the following. 

\begin{thm}\label{thm-counter}
There is a compact almost Hermitian $4$ manifold $(M,J,\omega)$ satisfying $\ell^{1,0}=\ell^{0,1}<\frac{b_1}{2}$, while $d\omega\neq0$.
\end{thm}


\begin{pf}
Choose $M=T^2\times S^2$. Define the almost complex structure on $M$ as follows:
\begin{itemize}
\item[(1)] On $TT^2$, set  $J\partial_x=\frac{1}{q}\partial_y$ and $J\partial_y=-q\partial_x$, where $q$ is a positive Morse function depending on $S^2$ and $(x,y)$ denotes the coordinate of $T^2$. 
\item[(2)] On $TS^2$, set $J|_{TS^2}=J_0$, where $J_0$ is the standard complex structure of $S^2$.   
\end{itemize}
Choose $\omega=p\cdot \omega_0+\omega_s$ as a nondegenerate $2$form. Here, $p$ is a non-constant positive function depending only on  $S^2$ and  $\omega_0,~\omega_1$ are the standard symplectic forms on $T^2,~S^2$ respectively. 
We can check that $J$ is $\omega$-compatible and $$d\omega\neq0.$$ 
Clearly,
$\theta:=dx+iq\cdot dy$ is a global $(1,0)$-form on $T^{*1,0}M$, whose dual is 
\[v:=\frac{1}{2}(\frac{\partial}{\partial x}-\frac{i}{q}\frac{\partial}{\partial y}).\]
Let $z$ a coordinate in a neighborhood $U\subset S^2$.
We write the structure equation
\begin{eqnarray}
(\bar\mu\theta)|_{T^2\times U}=\frac{1}{2q}\frac{\partial q}{\partial \bar z}\bar\theta\wedge d\bar z.\label{eqn-int-str}  
\end{eqnarray}
This shows that $\bar\mu \theta$ and $\mu\bar\theta$ are a.e. non-zero forms by the choice of Morse function.
Hence, a form $u\in\mathcal A^{0,1}$ can be locally expressed in $T^2\times U$ as
\[u=f_1\bar\theta+f_2 d\bar z,\]
where $f_1$ is a globally defined smooth function of $M$.
The equation $\mu(u)=0$ and \eqref{eqn-int-str}  guarantee that
\[\frac{\partial q}{\partial  z}\cdot f_1=0.\]
Since $q$ is a Morse function, $f_1\equiv0$. 
We locally expand the equation $\bar\partial u=0$ as 
\[\bar vf_2\cdot \bar\theta\wedge d \bar z=0.\]
That is $\bar v f_2=0$. Note that $v\bar v=\frac{1}{4}(\partial^2_x+\frac{1}{q^2}\partial^2_y)$ is an elliptic operator along $T^2$. The  maximum-principle asserts that $f_2$ is locally independent of $T^2$. Namely, $u$ is locally independent of $T^2$. For two  coordinates $z$ and $w$ in $U$ and $V$ respectively, one has the transformation $$f_2(x,y;z)=f_2(x,y,;w)\frac{\partial w}{\partial z},$$ in their intersection,  where the transition function $\frac{\partial w}{\partial z}$ depends only on $S^2$ and is non-zero. So, $u$ is the pull-back of  a $(0,1)$-form on $S^2$. 
Expand the equation $\bar\partial^*u=0$ as
\begin{eqnarray*}
0&=&\partial(f_2d\bar z\wedge \omega)\\
&=&\partial_zf_2\cdot p  (dz\wedge d\bar z\wedge\omega_0)+\partial_z p\cdot f_2(dz\wedge d\bar z\wedge\omega_0)\\
&=&\partial_z(f_2\cdot p)\cdot (dz\wedge d\bar z\wedge\omega_0).
\end{eqnarray*}
We get $\partial_z(f_2\cdot p)=0$ in $U$. By the coordinate change,
\[\partial_w(f_2(w)p(x,y;w))=\frac{\partial z}{\partial w}\partial_z(f_2(z)p(x,y;z))=0,\]
we rewrite as $$\partial_{S^2}(pu)=0,$$
where $\partial_{S^2}$ is the standard $\partial$-operator on $S^2$. 
We get $p \cdot u=0$ by $H^{1,0}(S^2)=0$. Since $p>0$, we have $u\equiv0$. Analogously, we also have that any $(1,0)$-form in $\ker(\bar\mu)\cap \ker(\bar\partial)$ is zero. That is, 
$h^{1,0}=\ell^{1,0}=0$ by Lemma \ref{lemma-h-ell}.
\end{pf}


\noindent
The following example shows that the even condition on the Betti numbers is not sufficient to make compact symplectic manifolds satisfy the Hard Lefschetz condition. 

\def\mfg{\mathfrak{g}}
\begin{exa}\label{exa-non-hlc}
Let $N$ be a nil-potent simple connected $4$ dimensional  Lie group, whose non-trivial structure equations on its associated Lie algebra $\mfg=\mbox{span}(e_1,e_2,e_3,e_4)$ are given by 
\[[e_1,e_2]=e_3,~[e_1,e_3]=e_4.\]
Then, on the dual algebra $\mfg^*$, it holds that
\[d e^1=0 ,  ~  d e^2=0,   ~  d e^3= - e^{1}\wedge e^2,    ~  d e^4= - e^{1}\wedge e^3.\]
Recall that the Lie group $N$ admits a lattice $ \Gamma$. Let $M = \Gamma\backslash N$ and
the form $\omega = e^{1}\wedge e^4 + e^{2}\wedge e^3$
is symplectic. It is clear that $b^1=2$. But $(M,\omega)$ does not satisfy the Hard Lefschetz condition, since 
$$[\omega\wedge e^1]= [ e^{1}\wedge e^2\wedge e^3 ] = [d (e^4\wedge e^{2})] = 0.$$
Moreover, for any compatible almost complex structure, we also have $ {h}^{1,0}=0$ by the inequalities $0\leq 2  h^{1,0}<b^1$. 
\end{exa}

\begin{pf}\textbf{( of Theorem \ref{lemma-b1-almost-kahler})}
Recall that $\Delta_d$ is a real operator, and each real-valued $1$-form in $\mathcal{A}^1_{\mathbb R}$ can be written as $u+\bar u$ for a unique element $u\in\mathcal{A}^{1,0}$, where $\mathcal{A}^1_{\mathbb R}$ is the space of real-valued $1$-forms. It suffices to show that each element $u+\bar u\in\ker(d)\cap\ker(d^*)\cap\mathcal{A}^1_{\mathbb R}$ corresponds to a unique element $u\in H^{1,0}_{Dol}$ on $(X,J,g,\omega)$.  

Combining  the condition $d^*(u+\bar u)=0 ,$   
\mbox{with the K\"ahler identities  }$[\Lambda,\partial]=i\bar\partial^*\mbox{ and }[\Lambda,\bar\partial]=-i\partial^*$, one obtains  $\Lambda(\bar\partial u-\partial\bar u)=0$. 
Together with the equation  $(\bar\partial u+\partial \bar u)=0$, we get  
\begin{equation}
\Lambda\bar\partial u=\Lambda\partial \bar u=0,\label{eqn-asd}
\end{equation}
i.e., $\partial\bar u$ is an imaginary-valued anti-self-dual $2$-form. 
We have the formula   
\begin{eqnarray*}       
\int_X\partial u\wedge\bar\partial\bar u&=&        -\int_X\bar\partial\partial u\wedge\bar u\\        &=&-\int_X(\partial^2\bar u+\mu\bar\partial\bar u)\wedge \bar u\\       
&=&\int_X\partial\bar u\wedge\partial\bar u+\int_X\bar\partial\bar u\wedge\mu\bar u\\        &=&\int_X\partial\bar u\wedge\partial\bar u-\int_X\bar\partial\bar u\wedge \partial u,   
\end{eqnarray*}   
where we use  the formula $\bar\partial\partial u=-\partial\bar\partial u-\bar\mu\mu u-\mu\bar\mu u=    \partial^2\bar u+\mu\bar\partial\bar u$ for the second equality and the formula $\partial u+\mu\bar u=0$ for the last one.   
The above formula and the anti-self-duality of $\bar\partial u$ yield 	\begin{equation}
\|\bar\partial u\|^2_{L^2}=\int_X\partial\bar u\wedge*\bar\partial u=-\int_X\bar\partial u\wedge\partial\bar u=\int_X\partial\bar u\wedge\partial\bar u=2\|\partial u\|^2_{L^2}.\label{eqn-parital-d}
\end{equation}  
On the other hand, we compute  
\begin{eqnarray*}     
\int_X(\Delta_du, u)&=&\int_X(du, du) +(d^*u,d^*u) \\
&=&\int_X(\bar\partial u,\bar\partial u)+(\partial u,\partial u)+(\bar\mu u,\bar\mu u)  = 4\int_X(\Delta_{\bar\mu} u,u), 
\end{eqnarray*}    
where we use the formulas $ d^*u=\partial^*u=0$, \eqref{eqn-asd} for the second equality and  the formulas $\bar\partial \bar u+\bar\mu u=0$, \eqref{eqn-parital-d} for the last one. Therefore, by the hypothesis of $\lambda_1$, it holds that $u\in\ker(d)\cap\mathcal{A}^{1,0}$, i.e., $u\in  H^{1,0}_{Dol}$.
\end{pf}

We end this paper by showing the following proposition for higher dimensional symplectic satisfying the Hard Lefschetz condition. For more results on the rank of Nijenhuis tensor on  almost complex manifolds, one can refer  \cite{CGGH21}, \cite{CPS}, \cite{Mu86}and \cite{ST23}. 
First,  we give a technique lemma. 
\begin{lem}\label{lemma-trivial}
Let $(M,J)$ be a compact almost complex $2n$-manifold  $n\geq3$. If the Nijenhuis tensor takes the maximal rank at some point of $M$, then $H^{1,0}_{Dol}=\{0\}$. 
\end{lem}
\begin{pf}
An easy   calculus yields $$rank_{\mathbb C}\mathcal{A}^{1,0}=rank_{\mathbb C}\mathcal{A}^{0,1}=n,~rank_{\mathbb C}\mathcal{A}^{2,0}=rank_{\mathbb C}\mathcal{A}^{0,2}=\frac{n(n-1)}{2}.$$ 
For an Hermitian metric $g$, it is clear that $\bar\partial^* u=0$ and $u\in\ker(\bar\partial)$ is equivalent to $u\in\ker(\Delta_{\bar{\partial}})$ for any $u\in \mathcal{A}^{1,0}$. 
Hence, all elements in $\tilde{H}^{1,0}$ are $\Delta_{\bar{\partial}}$-harmonic. Since the maximal rank condition is open, there exists a small open subset $U$ such that the restriction $N|_U$ is also of the maximal rank. 
When $n\geq3$, it is clear that $\bar\mu|_U:\mathcal{A}^{1,0}(U)\to\mathcal{A}^{0,2}(U)$ is either  isomorphic($n=3$) or injective($n>3$) by the algebraic calculus. The equation $\bar\mu(u)=0$ gives that $u|_U=0$ for any $u\in\mathcal{A}^{1,0}$.
Hence,    the unique-continuity-property for harmonic forms shows that $u=0$ for any $u\in{H}^{1,0}_{Dol}$. 

\end{pf}

\begin{prop}\label{prop-main-2} Let $(M,\omega,J)$ be a compact almost K\"ahler  manifold with dimension $\geq6$. 
If $b_1=2h^{1,0}\neq0$ for some compactible almost complex structure, then either $(M,\omega,J)$  is K\"ahlerian or the associated Nijenhuis tensor can not take the maximal rank at any point. 
\end{prop}
\begin{pf}
We argue by contradiction. Assume that $(M,\omega, J)$ is non-K\"ahlerian. 

If the associated Nijenhuis tensor $N$ takes maximal rank at a point $p$, then there exists a small neighborhood $U_p$ such that $N\mid_{U_p}$ is of maximal rank by the open condition of the maximal rank. 
Again, by the relation  $\bar\mu|_{U_p}:\mathcal{A}^{1,0}(U_p)\to\mathcal{A}^{0,2}(U_p)$ is either  isomorphic($n=3$) or injective($n>3$), any $u\in {H}^{1,0}_{Dol}$ is zero inside $U_p$. 
Thus, $u\equiv0$ by the unique continuity property, which contradicts to the assumption $b_1\neq0$. 
\end{pf}

College of Mathematics and Statistics, Chongqing University,
Huxi Campus, Chongqing, 401331, P. R. China

Chongqing Key Laboratory of Analytic Mathematics and Applications, Chongqing University, Huxi Campus, Chongqing, 401331, P. R.
China

E-mail: lindexie@126.com

\end{document}